\documentclass[12pt,reqno]{amsart}
\usepackage{microtype}

\usepackage[shortlabels]{enumitem}
\setlist[enumerate]{label={(\arabic*)}}

\usepackage{amssymb}
\usepackage[dvipsnames]{xcolor}
\usepackage
[colorlinks=true,linkcolor=Maroon,citecolor=OliveGreen,backref]
{hyperref}
\usepackage[abbrev,shortalphabetic]{amsrefs}  
\usepackage[capitalize]{cleveref} 
\crefname{equation}{}{}
\usepackage{tikz-cd}
\usepackage{tikz}

\usepackage[norefs, nocites]{refcheck}

\makeatletter
\newcommand{\refcheckize}[1]{%
  \expandafter\let\csname @@\string#1\endcsname#1%
  \expandafter\DeclareRobustCommand\csname relax\string#1\endcsname[1]{%
    \csname @@\string#1\endcsname{##1}\wrtusdrf{##1}}%
  \expandafter\let\expandafter#1\csname relax\string#1\endcsname
}
\makeatother

\refcheckize{\cref}
\refcheckize{\Cref}

\newtheorem{lemma}{Lemma}
\newtheorem{theorem}[lemma]{Theorem}
\newtheorem{proposition}[lemma]{Proposition}
\newtheorem{corollary}[lemma]{Corollary}

\theoremstyle{definition}
\newtheorem{remark}[lemma]{Remark}

\newcommand\opr[1]{\operatorname{#1}}

\renewcommand{\bar}{\overline}

\def\Sym{\opr{Sym}}

\newcommand\br[1]{{\left(#1\right)}}
\newcommand\floor[1]{\left\lfloor{#1}\right\rfloor}

\newcommand\pfrac[2]{\br{\frac{#1}{#2}}}
\newcommand{\gen}[1]{\left\langle{#1}\right\rangle}

\def\P{\opr{Prob}}

\usepackage{todonotes}

\begin{document}

\title{Dixon's asymptotic without CFSG}

\author{Sean Eberhard}
\address{Sean Eberhard, Mathematical Sciences Research Centre, Queen's University Belfast, Belfast BT7~1NN, UK}
\email{s.eberhard@qub.ac.uk}

\thanks{SE is supported by the Royal Society.}

\begin{abstract}
Without using the classification of finite simple groups, we show that the probability that two random elements of $S_n$ generate a primitive group smaller than $A_n$ is at most $\exp(-c(n \log n)^{1/2})$.
As a corollary we get Dixon's asymptotic expansion
\[
    1 - 1/n - 1/n^2 - 4/n^3 - 23/n^4 - \cdots
\]
for the probability that two random elements of $S_n$ (or $A_n$) generate a subgroup containing $A_n$.
\end{abstract}
\maketitle


\section{introduction}

We give a CFSG-free proof of the following result.

\begin{theorem}
    \label{thm:1}
    Let $G$ be the subgroup of $S_n$ generated by two random elements.
    The probability that $G$ is contained in a primitive subgroup of $S_n$ smaller than $A_n$ is bounded by $\exp(-c (n \log n)^{1/2})$ for some $c>0$.
\end{theorem}

This improves \cite{eberhard-virchow}*{Theorems~1.3 and 1.6}.
By combining with the results of \cite{dixon-asymptotics} we have the following corollary. (See also \cite{oeis}*{A113869}.)

\begin{corollary}
    The probability that two random elements of $A_n$ generate the group is
    \[
        1 - 1/n - 1/n^2 - 4/n^3 - 23/n^4 - 171/n^5 - \cdots.
    \]
    The same asymptotic expansion is valid for the probability that two random elements of $S_n$ generate at least $A_n$.
\end{corollary}

\section{Satisfaction probability for unimodal words}

Let $F_2 = F\{x,y\}$ be the free group on two letters $x, y$. We write $\{x,y\}^*$ for the set of positive words, i.e., the submonoid generated by $\{x,y\}$. Let $G = S_n = \Sym(\Omega)$ for $\Omega = \{1,\dots,n\}$.

\begin{proposition}
    Let $u, v \in \{x,y\}^*$ be distinct and let $w = uv^{-1} \in F_2$.
    Let $\ell = \ell(w) = \ell(u) + \ell(v)$ be the length of $w$.
    For a random evaluation $\bar w = w(\bar x,\bar y)$ with $\bar x, \bar y \in S_n$ uniformly random and independent, we have
    \[
        \P(\bar w = 1) \le (2\ell/n)^{\floor{n/2\ell}}.
    \]
\end{proposition}
\begin{proof}
    Write $w = w_1 \cdots w_\ell$ with $\ell > 0$ and $w_i \in \{x^{\pm1}, y^{\pm1}\}$ for each $i$.
    We may assume this expression is cyclically reduced.

    We use the query model for random permutations (see \cite{broder-shamir} or \cite{EJ}*{Section~A.1}).
    We gradually expose a random permutation $\pi \in \Sym(\Omega)$ by querying values of our choice.
    At every stage $\bar x$ and $\bar y$ are partially defined permutations.
    We may query the value of any $\pi \in \{\bar x^{\pm 1}, \bar y^{\pm 1}\}$ at any point $\omega \in \Omega$.
    If $\omega$ is already in the known domain of $\pi$, the known value is returned; this is a \emph{forced choice}.
    Otherwise, a random value is chosen uniformly from the remaining possibilities (the complement of the known domain of $\pi^{-1}$);
    this is a \emph{free choice}.
    If the result of a free choice is a point in the known domain of any of $\bar x^{\pm 1}, \bar y^{\pm 1}$ we say there was a \emph{coincidence}.
    It is standard and easy to see that this process results in uniformly random permutations $\bar x$ and $\bar y$ once all values are revealed.

    Begin by choosing any $\omega_1 \in \Omega$ and exposing the trajectory
    \[
        \omega_1^{\bar w_1}, \omega_1^{\bar w_1 \bar w_2}, \dots, \omega_1^{\bar w_1\cdots \bar w_\ell}.
    \]
    Let $E_1$ be the event that $\omega_1^{\bar w_1 \cdots \bar w_\ell} = \omega_1$.
    For this event to occur we claim it is necessary there was some coincidence among our queries of the form $\omega^{\bar w_t} = \omega_1$ (this is the crucial part of the argument).
    If $\ell(u) = 0$ or $\ell(v) = 0$ the argument is easy, so assume $u$ and $v$ have positive length.
    We may assume $w_1 = x$ and $w_\ell = y^{-1}$ since $w$ is cyclically reduced.
    If there is no coincidence of the given form, the trajectory of $\omega_1$ under $\bar u$ does not return to $\omega_1$, so $\omega_1$ cannot be added to the known domain of $\bar y$.
    Subsequently, during the negative part of the trajectory, unless there is a coincidence of the given form, $\omega_1$ can be added to the known domains of $\bar x^{-1}$ and $\bar y^{-1}$ only.
    Therefore at the final step $\omega_1$ is not in the known domain of $\bar y$,
    so if the final step is forced then the result is not $\omega_1$,
    and if the final step is free then the result is not $\omega_1$ by hypothesis.
    This proves the claim.

    Since the probability that any given free choice results in $\omega_1$ is at most $1 / (n - \ell)$, it follows by a union bound that
    \[
        \P(\omega_1^{\bar w} = \omega_1) \le \ell / (n-\ell).
    \]

    Conditional on the event $E_1$ choose a new point $\omega_2$ outside the trajectory of $\omega_1$, examine the trajectory of $\omega_2$, and so on.
    In general, at iteration $i$, conditional on the event $\bigcap_{j < i} E_j$ where $E_j = \{\omega_j^{\bar w} = \omega_j\}$, choose a point $\omega_i \in \Omega$ outside the union of the trajectories of $\omega_1, \dots, \omega_{i-1}$
    and query the trajectory of $\omega_i$.
    In order for the event $E_i = \{\omega_i^{\bar w} = \omega_i\}$ to occur it is necessary that there be a coincidence of the form $\omega^{\bar w_t} = \omega_i$.
    Therefore
    \[
        \P(\omega_i^w = \omega_i \mid E_1, \dots, E_{i-1}) \le \ell / (n - i \ell).
    \]

    Let $k = \floor{n/2\ell}$.
    Since the event $\{\bar w=1\}$ is contained in $E_1 \cap \cdots \cap E_k$, it follows that
    \[
        \P(\bar w = 1) \le \prod_{i = 1}^k \frac{\ell}{n-i\ell} \le \pfrac{2\ell}{n}^{\floor{n/2\ell}}.\qedhere
    \]
\end{proof}

\begin{remark}
    The proof above is essentially that of \cite{GHSSV}*{Section~3}.
    An error in that argument was identified in \cite{eberhard2017trivial},
    but the problem does not arise for words of the special form $w = uv^{-1}$,
    as explained in the third paragraph of the proof.
\end{remark}

\section{The order of the group}

Now let $\bar x, \bar y \in S_n$ be uniformly random and let $G = \gen{\bar x,\bar y}$.

\begin{proposition}
    There is a constant $c > 0$ such that
    \[
        \P(|G| < \exp(c(n \log n)^{1/2})) \le \exp(-c(n \log n)^{1/2}).
    \]
\end{proposition}
\begin{proof}
    Consider the elements of $G$ of the form $\bar u$ with $u \in \{x,y\}^*$ and $\ell(u) < r$ (for some $r$).
    The number of such $u$ is $1 + 2 + \cdots + 2^{r-1} = 2^r-1$.
    Applying the previous proposition, the probability that any two such $\bar u$ and $\bar u'$ are equal is bounded by
    \[
        4^r (4r/n)^{\floor{n/4r}} \le \exp(c_1 r - c_2 (n/r) \log(n/r)))
    \]
    for some constants $c_1, c_2 > 0$.
    Choosing $r = c_3 (n \log n)^{1/2}$ for a small enough constant $c_3>0$, we obtain a bound of the required form.
    Failing this event, $|G| \ge 2^r - 1$, so the result is proved.
\end{proof}

A beautiful recent result of Sun and Wilmes~\cite{sun--wilmes--abstract, sun--wilmes} (building on seminal work of Babai~\cite{babai--annals}) classifies primitive coherent configurations with more than $\exp( C n^{1/3} (\log n)^{7/3})$ automorphisms.
A corollary is a CFSG-free determination of the uniprimitive subgroups of $S_n$ of order greater than the same bound.
Much stronger bounds for the order of $2$-transitive groups have been known for a long time~\cite{babai-2-trans,pyber-2-trans}.
Thus we know there are at most two conjugacy classes of primitive maximal subgroups $M < S_n$ apart from $A_n$ such that $|M| > \exp(C n^{1/3} (\log n)^{7/3})$,
and each satisfies $|M| = \exp \br{O(n^{1/2} \log n)}$.
Since the number of pairs of permutations lying in a common conjugate of a maximal subgroup $M$ is at most $1/[S_n:M]$,
\Cref{thm:1} follows.

\begin{remark}
    This proof was essentially anticipated in \cite{babai-JCTA}*{Remark~1}.
\end{remark}

\bibliography{refs}
\end{document}